\documentclass[british,smallextended]{svjour3}
\usepackage[T1]{fontenc}
\usepackage[latin9]{inputenc}
\usepackage{array}
\usepackage{textcomp}
\usepackage{url}
\usepackage{multirow}
\usepackage{amsmath}
\usepackage{graphicx}
\PassOptionsToPackage{version=3}{mhchem}
\usepackage{mhchem}
\usepackage[authoryear]{natbib}

\makeatletter

\providecommand{\tabularnewline}{\\}

\usepackage{changebar}
\providecommand{\lyxadded}[3]{}

\renewcommand{\lyxadded}[3]{%
  {\protect\cbstart\color{lyxadded}{}#3\protect\cbend}%
}

\RequirePackage{fix-cm}

\smartqed  

\renewcommand\[{\begin{equation}}
\renewcommand\]{\end{equation}} 

\makeatother

\usepackage{babel}

\begin{document}
\title{Embedding Reservoirs in Industrial Models to Exploit their Flexibility\thanks{This research was supported by the public service of Wallonia, within
the framework of the InduStore project (grant 1450300).}}
\author{Thibaut Cuvelier}
\institute{T. Cuvelier \at  L2S (\emph{laboratoire des signaux et systèmes}),
CentraleSupélec, Gif-sur-Yvette, France\\
This work was carried out while the author was a research engineer
at the university of Liège, Belgium.\\
\email{thibaut@tcuvelier.be}}
\date{Received: date / Accepted: date}
\maketitle
\begin{abstract}
In the context of energy transition, industrial plants that heavily
rely on electricity face more and more price volatility. To continue
operating in these conditions, the directors become continually more
willing to increase their flexibility, i.e. their ability to react
to price fluctuations. This work proposes an intuitive methodology
to mathematically model electro-intensive processes in order to assess
their flexibility potential. To this end, we introduce the notion
of \emph{reservoir}, a storage of either material or energy, that
allows models based on this paradigm to have interpretations close
to the physics of the processes. The design of the reservoir methodology
has three distinct goals: ($i$) to be easy and quick to build by
an energy-sector consultant; ($ii$) to be effortlessly converted
into mixed-integer linear or nonlinear programs; ($iii$) to be straightforward
to understand by nontechnical people, thanks to their graphic nature.
We apply this methodology to two industrial case studies, namely an
induction furnace (linear model) and an industrial cooling installation
(nonlinear model), where we can achieve significant cost savings.
In both cases, the models can be quickly written using our method
and solved by appropriate solver technologies. 
\end{abstract}

\section{Introduction\label{sec:Introduction}}

Current industries tend to consume large quantities of energy, often
electricity, during several stages of their processes: in aluminium
production, to extract the material from alumina; in cement making,
to crush the limestone and heat the kiln; in paper fabrication, to
pulp the ground wood. An industrial site is said to be electro-intensive
when the costs of energy compose a large part of the retail price.
In this case, the dependency on electricity prices is very high in
order to remain competitive. The literature mostly focuses on electricity
in this context, as its price can be very volatile, as opposed to
other fuels like natural gas.

Nevertheless, these industrial sites can often tune their processes
in order to decrease their electricity consumption during the most
expensive periods; they could also benefit from selling flexibility
services to the electrical network, by reducing their consumption
when the operator asks for it~\citep{Asadinejad2017}. To achieve
this cost reduction, they can make use of decision-support systems
based on mathematical modelling of their processes. For example, the
heating of a cement kiln can be reduced when the price is too high,
or completely stopped if the revenue from the cement is not sufficient
to cover the production costs. However, the achieved consumption drop
is frequently the result of a trade-off:
\begin{itemize}
\item in some cases, the system is operated differently to use another fuel
(\emph{fuel switching}), and this has little to no impact on the process.
For instance, a cement kiln can switch from natural gas to flare gas
or other industrial by-products~\citep{Bosoaga2009}; in ethylene
production, production plants can typically use several fuels including
naphtha or Diesel oil~\citep{Petracci1996,Han2020};
\item in others, the only flexibility is to stop the process and lose some
production (\emph{load shedding}), which means that some orders cannot
be fulfilled. For example, in aluminium smelting, any reduction of
the current fed into the electrolysis potlines implies a loss of production~\citep{Todd2008};
\item other flexibility levers can be available, such as \emph{load shifting}
or \emph{scheduling}, depending on the exact process, with varying
degrees of impact on the production and the consumption. The paper
industry falls in this category, with production plan reordering~\citep{Dejemeppe2016}.
\end{itemize}
Indeed, the plant directors would welcome global models, describing
the whole set of processes, that allow them to make the best decisions
for their energy consumption based on price forecast. To this end,
each industrial process in the factory should be modelled to provide
a consumption-production mathematical formulation. These models do
not need to offer highly detailed insights into the machines: the
goal is not to have a real-time control of the physics thereof. In
fact, our aim is rather to get an \emph{estimation} of the consumption
of the process depending on how it is operated. Moreover, the staff
might not have the required knowledge of operational research to build
such models: they have to rely on external consultants to do so.

This article proposes a generic paradigm to help conceiving such approximate
formulations, so that they can be used to characterise the flexibility
of a given process by the means of mathematical optimisation. The
low-level basic block of this paradigm is the \emph{reservoir}, which
yields simple optimisation models while having a great expressive
power.

For example, an oven might be modelled as two such reservoirs: a material
tank (expressed in tonnes), and an energy storage (in joules). Over
time, the oven loses some thermal energy, which can be modelled as
losses from the energy reservoir. The oven is mainly controlled by
its temperature, constrained to some operating range; this temperature
is determined as a function of both the material and the energy reservoirs.
The expressive power of the proposed paradigm comes from the fact
that a reservoir may represent an actual physical storage or be more
abstract (such as an energy storage or any quantity useful for modelling).

The complete oven model must also contain other blocks, named \emph{processes}:
these can impact the reservoirs, such as heating the oven. Only those
processes might consume energy, such as electricity or natural gas
for oven burners. A set of reservoirs only represents the state of
the industrial machine, while processes interact with this state.

Using reservoir-based models for industrial processes has several
advantages over typical ad-hoc models, especially in the context of
building many such models in a short amount of time:
\begin{itemize}
\item Reservoir-based models are built from a small number of intuitive
components (described in Section~\ref{sec:Reservoir-taxonomy}).
Consultants could quickly grasp the main ideas.
\item Reasoning about such a model is easier than with pure mathematical
notations, as this paradigm suggests representing the processes as
legible diagrams. This graphical representation helps communicate
with technical experts of the industrial process and nontechnical
managers.
\item The resulting optimisation models are often linear mixed-integer (MILPs).
\item Conversion into computer code is straightforward once the building
blocks are implemented. A software implementation may even propose
a graphical interface to build the complete optimisation model, similar
to MATLAB Simulink~\citep{Simulink2020}. Almost no adaptations are
then required to fit the parameters of such a model to known data
(as highlighted in Section~\ref{subsec:Results-Induction-furnace}).
\end{itemize}
This article first details the developed reservoir taxonomy in Section~\ref{sec:Reservoir-taxonomy},
and some usage examples in Section~\ref{sec:Reservoir-models}, including
diagrams. Even though the resulting problems are generally linear,
more complex (i.e. nonlinear) behaviours can be implemented, as explained
in Section~\ref{sec:Nonlinearity}. Based on this experience, we
build a typology of processes in Section~\ref{sec:Process-typology}.
Our numerical results are presented in Section~\ref{sec:Results},
both for fitting a reservoir model to industrial data and for assessing
the flexibility within an industrial site. We conclude in Section~\ref{sec:Conclusion}.

\section{Reservoir taxonomy \label{sec:Reservoir-taxonomy}}

To develop the aforementioned simplified models, four kinds of building
blocks are needed. Each of them is directly associated to a mathematical
formulation. In the following notations, \textbf{boldface} indicates
optimisation variables, as opposed to constants, typeset in roman.
\begin{itemize}
\item The \emph{reservoir}, which is as close to the intuition of storage
as possible: the level $\mathbf{s}_{t}$ is only impacted by the inflow
$\mathbf{in}_{t}$ and the outflow $\mathbf{out}_{t}$. 
\[
\mathbf{s}_{t+1}=\mathbf{s}_{t}+\mathbf{in}_{t}-\mathbf{out}_{t}.
\]
These inflow and outflow variables link the reservoir to other parts
of a global plant model. In the paper industry, for instance, the
excess of pulp production may be stored before it is used by the paper
machines, which indicates some decoupling between the production and
the consumption \textemdash{} in other words, a source of flexibility.
\item The \emph{decaying reservoir}, whose distinct feature is to have leaks
$\mathbf{decay}_{t}$ in its content. 
\[
\mathbf{s}_{t+1}=\mathbf{s}_{t}+\mathbf{in}_{t}-\mathbf{out}_{t}-\mathbf{decay}_{t}.
\]
Those leaks can be modelled with any kind of mathematical relationship:
they can depend on the state of the reservoir $\mathbf{s}_{t}$, but
also on the state of \emph{other }reservoirs. For example, an industrial
oven naturally loses some thermal energy over time.
\item The \emph{observer}, which is not a reservoir per se. Its state is
a function of other reservoirs' state, denoted by the variables $\mathbf{u}_{t},\mathbf{v}_{t}\dots$:
\[
\mathbf{s}_{t}=f\!\left(\mathbf{u}_{t},\mathbf{v}_{t}\dots\right).
\]
It is mainly useful to implement bounds based on the state of a reservoir,
such as the temperature in an oven, when it is modelled with both
an energy and a mass reservoirs.
\item On top of these, \emph{external processes} must be added to impact
the state of the reservoirs, by imposing some value to the inflow
and outflow variables. Two examples are heaters (they increase the
level of a given heat reservoir) and chemical reactions (they transform
some products into some others).
\end{itemize}
Each of these reservoirs may impose constraints on its level variable
$\mathbf{s}_{t}$, such as bounds (like minimum and maximum temperatures),
ramping constraints (e.g., to limit the temperature variations), or
process-dependent constraints (an oven cannot be tapped before its
content is melted, for instance).

The flows between reservoirs may need to be coupled, especially when
a given process is modelled as multiple reservoirs. This kind of constraint
is common for heat transfer, for example: the quantity of heat transferred
is proportional to the flow of heat-transfer fluid.

Aggregating those building blocks constructs models for the whole
process. The final step to mathematically formulate the complete plant
is to assemble the various systems in one model.

\section{Reservoir models \label{sec:Reservoir-models}}

This section shows a series of models that can be obtained with the
methodology described in Section~\ref{sec:Reservoir-taxonomy}. All
the examples are taken from the industry, with HVAC (Section~\ref{subsec:HVAC-model})
being present in many sites; ovens (Section~\ref{subsec:Oven-model})
are mostly present in metallurgy, whereas electrolysis (Section~\ref{subsec:Electrolysis-model})
is a major process in chemistry. These examples are all electro-intensive
processes: a typical metallurgical EAF is around 30 MW, while it is
not uncommon to have electrolysis potlines of 50 MW.

Our goal is to produce \emph{simple and approximate} models: they
are used to detect the flexibility potential of a plant, not to perform
complex real-time control on the processes. This is why linear models
are usually good enough for our purposes. Would they fail in giving
a good enough estimation, nonlinearities would have to be introduced
(as done in Section~\ref{sec:Nonlinearity}).

One of the main advantages of the proposed methodology is that the
models can be summarised by drawings, while retaining many details.
This allows easy designing of models, while the actual formulation
for the links between the reservoirs must be dealt with in more details
outside the drawings.

\subsection{HVAC \label{subsec:HVAC-model}}

Heating, ventilation, and air conditioning (HVAC) is a kind of system
often found in the industry (food processing, supermarket warehouses,
pharmaceutical plants, etc.), but also in offices. This kind of process
is usually not the most electro-intensive that can be found, but the
industrial partners often agree with its flexibilisation (as opposed
to their main business).

\begin{figure}
\begin{centering}
\includegraphics[scale=0.4]{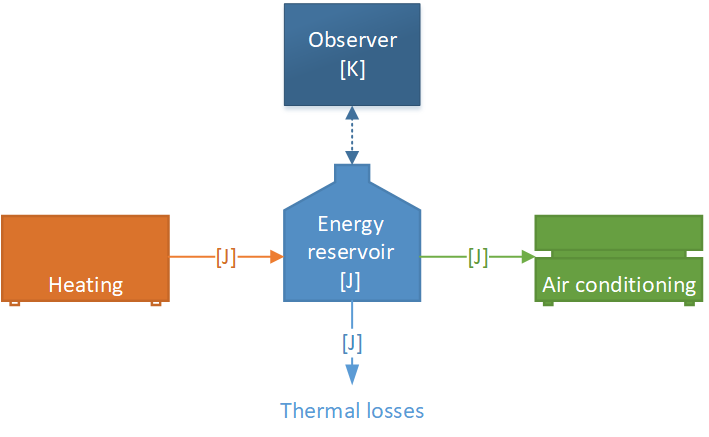}
\par\end{centering}
\caption{\label{fig:Reservoirs-HVAC}An HVAC installation modelled with reservoirs.}
\end{figure}

The proposed model works on a single \emph{energy reservoir}, corresponding
to the room whose temperature is controlled (or a set of rooms, or
an industrial shed), as shown in Figure~\ref{fig:Reservoirs-HVAC}.
Exploiting the thermal inertia of the area, the HVAC system can be
turned off to save on energy while keeping an acceptable temperature.

A series of \emph{external processes} heats or cools down the room
by adding or removing energy (heating and air conditioning): they
are the inputs and outputs of the reservoir. Moreover, the reservoir
level decreases due to thermal losses, which depend on the temperature
difference between the inside and the outside. An \emph{observer reservoir}
is used to control directly the room temperature.

The corresponding mathematical formulation is the following.
\begin{itemize}
\item The standard \emph{decaying reservoir} model is used for the energy.
It defines the main decision variable of this system: the energy within
the room, denoted by $\mathbf{E}_{t}$ at time $t$. Its state equation
takes into account the heater ($\mathbf{heating}_{t}$), the chiller
($\mathbf{ac}_{t}$), and thermal losses ($\mathbf{losses}_{t}$;
these may very well be positive, if the outside temperature is higher
than that of the room). 
\begin{equation}
\mathbf{E}_{t+1}=\mathbf{E}_{t}+\mathbf{heating}_{t}-\mathbf{ac}_{t}-\mathbf{losses}_{t}.\label{eq:HVAC-model-state}
\end{equation}
The decay quantity $\mathbf{losses}_{t}$ can be approximated as directly
proportional to the temperature difference between the room (denoted
by $\mathbf{T}_{t}$) and the outside ($T_{t}^{\mathrm{out}}$) \citep{Crombie2006},
the coefficient being denoted by $k$. 
\[
\mathbf{losses}_{t}=k\times\left(\mathbf{T}_{t}-T_{t}^{\mathrm{out}}\right).
\]
\item Using an \emph{observer}, the temperature is linked with the energy
through the specific heat of the air $C_{p}$ and the mass of air
in the room $m$ (which is considered constant): 
\[
\mathbf{E}_{t}=C_{p}\times m\times\mathbf{T}_{t}.
\]
This temperature has strict bounds ($T_{\min}$ and $T_{\max}$) that
are imposed inside the observer: 
\begin{equation}
T_{\min}\leq\mathbf{T}_{t}\leq T_{\max}.\label{eq:HVAC-model-temperature-bounds}
\end{equation}

\begin{remark}
As the mass of air $m$ is considered fixed, those constraints are
linear. This approximation works well, because the volume of air to
heat or cool down sees almost no variations in typical HVAC conditions.
Section~\ref{subsec:Industrial-cooling} considers a use case where
this hypothesis no more holds.
\end{remark}
\item \emph{External processes} relate the heater and the chiller to their
energy consumption (for example, electricity, $\mathbf{electricity}_{t}$,
and natural gas, $\mathbf{gas}_{t}$; a superscript $h$ denotes the
chiller, while $ac$ is used for the air conditioning). The exact
relationships $f$ and $g$ are not made explicit here: 
\[
\mathbf{heating}_{t}=f\!\left(\mathbf{electricity}_{t}^{\mathrm{h}},\mathbf{gas}_{t}^{\mathrm{h}}\right),
\]
\[
\mathbf{ac}_{t}=g\!\left(\mathbf{electricity}_{t}^{\mathrm{ac}},\mathbf{gas}_{t}^{\mathrm{ac}}\right).
\]
These relationships may be either exact or approximated, depending
on the need for precision. A basic model for these relations could
use an efficiency or a coefficient of performance to relate the energy
consumption and the produced effect (either heating or air conditioning).
\item The initial conditions indicate the state of the room to model at
the beginning of the optimisation horizon, which gives the initial
energy ($\mathbf{E}_{0}$) based on the initial temperature ($T_{0}$):
\[
\mathbf{E}_{0}=C_{p}\times m\times T_{0}.
\]
\end{itemize}
In a practical application, to fit this kind of model on experimental
data, three important variables stand out: the energy that is pushed
into or extracted from the room, the measured temperature, and the
outside temperature. This data is sufficient to get the values for
all the needed parameters ($k$ for the losses and the product $C_{p}\,m$
for the temperature). With linear regressions, the parameter $C_{p}\,m$
is determined as the proportion coefficient between the energy and
the temperature variations. With the same technique, the loss coefficient
$k$ may be determined from~\eqref{eq:HVAC-model-state}.

This model is simpler than existing ones in the literature~\citep{Kusiak2010,Lee2015},
but similar to models used in the context of flexibility~\citep{Short2019}.
\begin{remark}
This formulation neglects the effects rooms may have on each other.
Those may be included in such a model by using an energy reservoir
per room. The losses would then depend on the connected rooms. \\
Also, the average specific heat variations in the room are neglected,
as well as the temperature uniformity, i.e. the air is considered
homogeneous.
\end{remark}

\subsubsection{Industrial cooling \label{subsec:Industrial-cooling}}

A very similar model can be built for industrial cooling. The main
differences are that no heating is performed, and that cooling can
be done with two processes: either a cooling tower (whose efficiency
depends on the outside temperature) or a chiller.

\begin{figure}
\begin{centering}
\includegraphics[scale=0.4]{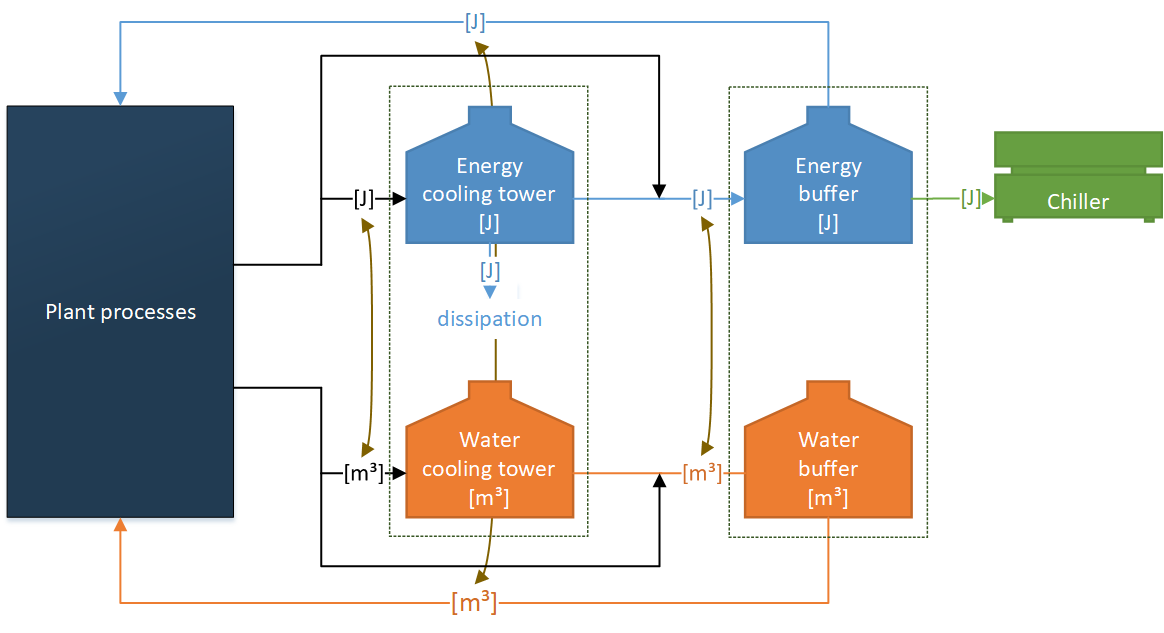}
\par\end{centering}
\caption{\label{fig:Reservoirs-Cooling}Industrial cooling modelled with reservoirs.}
\end{figure}

This model is more complicated due to the heat exchanges between the
three components (the processes to cool and the two cooling mechanisms):
the mass of heat- transfer fluid cannot be assumed to be constant.
Indeed, the flows of hot and cooled water can be variable, but also
the split between the two cooling processes. Each part of the model
has two reservoirs: one for heat and the other for water, as shown
in Figure~\ref{fig:Reservoirs-Cooling}. The flows between those
reservoirs must be coupled by the means of temperature and the water's
heat capacity $C$: 
\[
\mathbf{flow}_{\mathrm{heat}}=C\times\mathbf{T}\times\mathbf{flow}_{\mathrm{water}}.
\]

The model is thus the following. Subscript $CT$ denotes the cooling
tower, $p$ the processes that must be cooled, and $B$ the buffer
to which the chiller is connected. The arrow $\to$ denotes a flow
between two components.
\begin{itemize}
\item Neither the water nor the heat reservoirs are decaying, as this decay
is usually very small in these applications. 
\[
\mathbf{water_{CT,t+1}}=\mathbf{water_{CT,t}}+\mathbf{water_{p\to CT,t}}-\mathbf{water_{CT\to B,t}}.
\]
\[
\mathbf{water_{B,t+1}}=\mathbf{water_{B,t}}+\mathbf{water_{p\to B,t}}-\mathbf{water_{CT\to B,t}}.
\]
\begin{align}
\mathbf{E}_{CT,t+1} & =\mathbf{E}_{\mathrm{CT},t}+C_{\mathrm{water}}\times T_{p,t}\times\mathbf{water_{p\to CT,t}}\nonumber \\
 & -\mathbf{E}_{\mathrm{CT\to B},t}-f\!\left(\mathbf{P}_{\mathrm{CT},t}\right).\label{eq:IC-model-budget-energy-CT}
\end{align}
\begin{align}
\mathbf{E}_{B,t+1} & =\mathbf{E}_{B,t}+C_{\mathrm{water}}\times T_{\mathrm{p},t}\times\mathbf{water_{p\to B,t}}+\mathbf{E}_{\mathrm{CT\to B},t}\nonumber \\
 & -C_{\mathrm{water}}\times\mathbf{T}_{\mathrm{B},t}\times\mathbf{water_{B\to p,t}}-g\!\left(\mathbf{P}_{\mathrm{chiller},t}\right).\label{eq:IC-model-budget-energy-chiller}
\end{align}
\item The total water and energy from the processes must correspond to the
input scenario. 
\[
\mathbf{water_{p\to CT,t}}+\mathbf{water_{p\to B,t}}=\mathrm{water}_{\mathrm{p},t}.
\]
\item The flows of water and energy are linked. 
\begin{equation}
\mathbf{E}_{\mathrm{CT},t}=C_{\mathrm{water}}\times\mathbf{water_{CT,t}}\times\mathbf{T}_{\mathrm{CT},t}.\label{eq:IC-model-concordance-CT}
\end{equation}
\begin{equation}
\mathbf{E}_{\mathrm{B},t}=C_{\mathrm{water}}\times\mathbf{water_{B,t}}\times\mathbf{T}_{\mathrm{B},t}.\label{eq:IC-model-concordance-B}
\end{equation}
\begin{equation}
\mathbf{E}_{\mathrm{CT\to B},t}=C_{\mathrm{water}}\times\mathbf{water_{CT\to B,t}}\times\mathbf{T}_{\mathrm{CT},t}.\label{eq:IC-model-concordance-flow}
\end{equation}
\item The buffer's temperature must be within acceptable bounds for the
processes, as its water is directly sent back to the processes. 
\[
T_{\min}\leq\mathbf{T}_{\mathrm{B},t}\leq T_{\max}.
\]
\end{itemize}
This model is simpler than existing ones~\citep{Soderman2010}, even
for similar applications, and has fewer parameters to fit~\citep{Peesel2019}.
Actually, in Section~\ref{subsec:Results-Industrial-cooling}, we
only use manufacturer-provided data to perform flexibility estimation.
The real performance may significantly differ from the specifications~\citep{Peesel2017},
but such a precision in the numerical results is usually not required
for our use case.

A highly similar model can be used for petrochemical chains, like
ethylene production~\citep{Han2020}. More complex models can still
be suited to mixed-integer linear formulations, but cannot be expressed
as reservoirs; they also come with higher computational costs with
traditional optimisation tools~\citep{Geng2020}. 

\subsection{Oven \label{subsec:Oven-model}}

\begin{figure}
\begin{centering}
\includegraphics[scale=0.4]{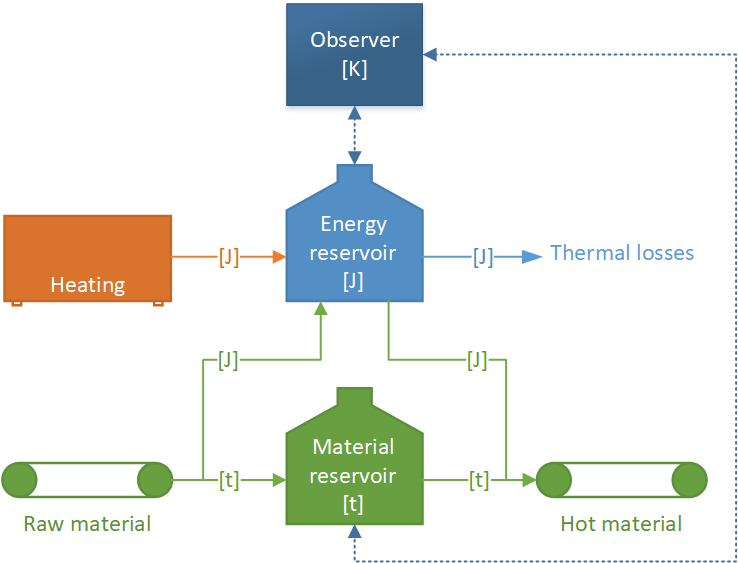}
\par\end{centering}
\caption{\label{fig:Reservoirs-oven}Industrial oven as reservoirs.}
\end{figure}

A more complex example of reservoir modelling is the industrial oven
(see Figure~\ref{fig:Reservoirs-oven}), which heats large quantities
of material (several tonnes) to high temperatures (around 1000 °C),
usually with a high thermal inertia (which allows for turning off
or reducing the heating from time to time). Different kinds of mechanisms
can be used to heat the materials: gas or oil burners (for relatively
low temperatures), electric arcs (for conductive materials), etc.

An oven model is similar to HVAC in that it has an energy reservoir.
However, the main difference is that the quantity of material to heat
may vary significantly, and thus cannot be neglected. This material
is modelled as a second reservoir, coupled with the first one.

The relationships in this model are more intricate than in the previous
ones: quantity of material, energy, and temperature are tightly and
nonlinearly intertwined. For instance, when heating, the impact on
temperature is not direct: for the same quantity of heating, the impact
on temperature is less when the oven contains a large quantity of
material than when it is almost empty.

In order for the oven to meet its operational goals, it must provide
the needed quantities of material at the right temperature. When material
is removed from the oven, it has lost an associated quantity of energy
\textemdash{} but not temperature. However, when new material is inserted
into the oven, this matter is at the outside temperature; it brings
some energy into the oven, but lowers the overall temperature.

All in all, the obtained mathematical formulation is the following.
\begin{itemize}
\item A first \emph{reservoir} considers the material, $\mathbf{m}_{t}$:
some may be added ($\Delta\mathbf{m}_{t}^{+}$) or removed ($\Delta\mathbf{m}_{t}^{-}$)
at any time step. 
\[
\mathbf{m}_{t+1}=\mathbf{m}_{t}+\Delta\mathbf{m}_{t}^{+}-\Delta\mathbf{m}_{t}^{-}.
\]
\item A \emph{decaying reservoir} is used for the energy. Adding or removing
material has an impact on the energy content (with the specific heat
$C_{p}$, as previously). The withdrawn matter is at the oven temperature;
however, the added material has a lower temperature, namely $T_{t}^{\mathrm{out}}$:
\[
\mathbf{E}_{0}=0,
\]
\begin{eqnarray}
\mathbf{E}_{t+1} & = & \mathbf{E}_{t}+\mathbf{heating}_{t}-\mathbf{losses}_{t}\label{eq:Oven-model-state-energy}\\
 & + & C_{p}\times\Delta\mathbf{m}_{t}^{+}\times T_{t}^{\mathrm{out}}\nonumber \\
 & - & C_{p}\times\Delta\mathbf{m}_{t}^{-}\times\mathbf{T}_{t}\nonumber 
\end{eqnarray}
The losses take the same form as previously, being linked to the
temperature difference with the exterior: 
\begin{equation}
\mathbf{losses}_{t}=k\times\left(\mathbf{T}_{t}-T_{t}^{\mathrm{out}}\right).\label{eq:Oven-model-losses}
\end{equation}
\item Likewise, an \emph{observer} gives the temperature, but this time
\emph{nonlinearly}: 
\begin{equation}
\mathbf{E}_{t}=C_{p}\times\mathbf{m}_{t}\times\mathbf{T}_{t}.\label{eq:Oven-model-temperature}
\end{equation}
This nonlinearity is not a problem for the temperature bounds, as
the constraint~\eqref{eq:HVAC-model-temperature-bounds} can be rewritten
in joules-kilograms \emph{without} using constraint~\ref{eq:Oven-model-temperature}:
\begin{equation}
C_{p}\times T_{\min}\times\mathbf{m}_{t}\leq\mathbf{E}_{t}\leq C_{p}\times T_{\max}\times\mathbf{m}_{t},\label{eq:Oven-model-temperature-bounds}
\end{equation}
\item The heater is still an \emph{external process} that consumes energy:
\[
\mathbf{heating}_{t}=f\!\left(\mathbf{electricity}_{t},\mathbf{gas}_{t}\right).
\]
\end{itemize}
As opposed to the HVAC model, this formulation cannot be linear due
to the losses, as they involve the temperature~\eqref{eq:Oven-model-losses},
i.e. the ratio between the energy and the mass. These issues are discussed
in Section~\ref{sec:Nonlinearity}.
\begin{remark}
This model exploits the hypothesis that the specific heat $C_{p}$
remains constant with temperature. Also, it does not consider phase
change, i.e. it only allows for heating material, not to melt it.
A large quantity of energy being required for the material to change
phase, adding this possibility in the model would require another
reservoir that specifically deals with the molten part.
\end{remark}

\subsubsection{Induction furnace \label{subsec:Induction-furnace}}

Induction furnaces are a kind of industrial oven. They work by magnetic
induction: the bucket containing the metal to melt (such as cast iron)
is surrounded by a a coil through which high-voltage alternative currents
are sent. The created magnetic field induces eddy currents in the
metal, which in turn heats it. The power dissipated in the metal is
directly proportional to the electrical power fed into the circuit
$\mathbf{P}_{t}^{\mathrm{elec}}$, with a constant ratio $\alpha$,
as shown in Section~\ref{sec:Appendix-IF}.

\begin{figure}
\begin{centering}
\includegraphics[scale=0.4]{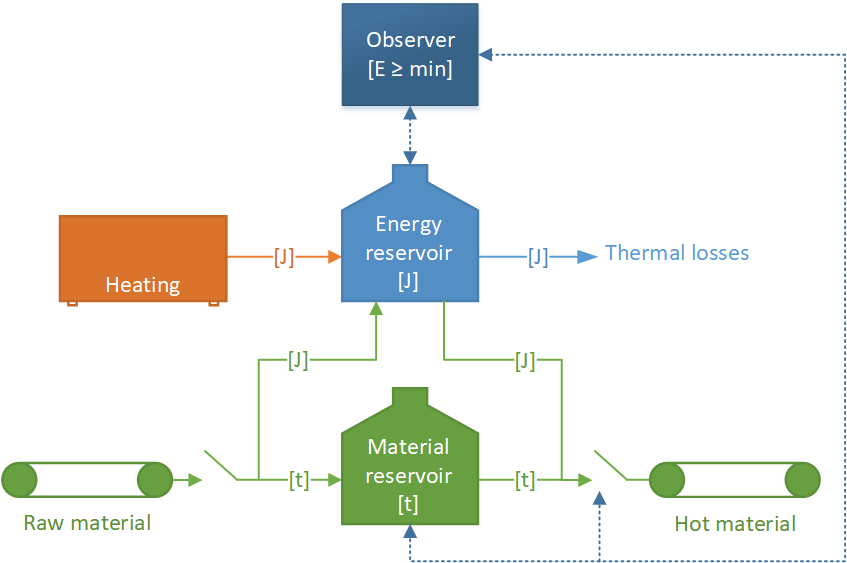}
\par\end{centering}
\caption{\label{fig:Reservoirs-IF}Induction furnaces modelled with reservoirs.}
\end{figure}

A reservoir model that fits this kind of furnace would be made up
of two reservoirs: a mass reservoir that evolves at discrete time
steps, and an energy decaying reservoir (filled and emptied at the
same time as the mass reservoir). There is no need for a temperature
observer, but only a ``binary'' observer that indicates whether
the required total energy has been transferred to the metal at the
end of the batch. The model is shown in Figure~\ref{fig:Reservoirs-IF}.

The corresponding model in mathematical form is the following:
\begin{itemize}
\item The mass reservoir is fixed for a given batch: 
\[
\mathbf{m}_{t}=m_{0}.
\]
\item The energy reservoir is filled by electrical currents: 
\[
\mathbf{E}_{t+1}=\mathbf{E}_{t}+\alpha\,\mathbf{P}_{t}^{\mathrm{elec}}.
\]
\item The observer reservoir imposes that, at the end of the batch, the
total energy in the bucket is at least at the required level: 
\[
\mathbf{E}_{T}\geq f\!\left(m_{0}\right).
\]
\item Initially, the energy reservoir indicates that the material is at
the exterior temperature:
\[
\mathbf{E}_{0}=C_{p}\times m_{0}\times T_{0}^{\mathrm{out}}.
\]
\end{itemize}
Typical models for such furnaces deal with the details of the electromagnetic
field required to melt the metal~\citep{Abubakre2009,Naar2013},
which is too cumbersome for our application.

\subsection{Electrolysis \label{subsec:Electrolysis-model}}

The chemical industry and metallurgy often use electrolysis, with
aluminium smelting being a prime example. A potline is typically made
of hundreds of individual baths, where the electrochemical reactions
take place. Each of them can be controlled independently and turned
off to save energy, but with a loss of production \citep{Todd2008}.

In front of the potline, there is a single transformer from AC to
DC current; its output has a relatively low voltage (around 5 V, usually),
but very high amperage (multiple hundreds of kiloamperes are not rare)~\citep{Todd2008,AlcoaInc.,BALCO}.
After the transformer, all pots resemble ovens (as in Section~\ref{subsec:Oven-model}):
they are heated by the current that flows through them, so that their
bath remains in a given temperature range where the electrochemical
reaction can happen.

The main difference with ovens is that the contents of the pots have
smaller fluctuations: inputs are continuous, so the mass increases
at a predictable and constant rate; the product is siphoned off periodically,
typically once per day~\citep{Todd2008,AlcoaInc.,BALCO}. However,
the contents of the bath continuously evolve: the alumina and the
carbon anode are transformed into aluminium and gas, the actual reactions
being $\text{Al}_{2}\text{O}_{3}\ce{ + 3 C \to\,2 Al + 3 CO}$ and
$\text{2 }\text{Al}_{2}\text{O}_{3}\ce{ + 3 C \to\,4 Al + 3 CO}_{2}$.
As such, $415\,\mathrm{kg}$ of carbon anode (not fed continuously,
the anode being replaced after a few weeks) are consumed to produce
$1\,\mathrm{t}$ of metallic aluminium, and this carbon is ejected
as exhaust gas (carbon mono- and dioxide)~\citep{TheNewZealandInstituteofChemistry2008}.
Nevertheless, these exhausts are neglected in our model.

The flexibility impacts the production in a slightly more complicated
way than in the aforementioned processes: it depends (approximately
linearly) on the average \emph{current} through the pot, as long as
the temperature is in the right range (which is controlled by voltage)~\citep{Molina-Garcia2011}.

\begin{figure}
\begin{centering}
\includegraphics[scale=0.4]{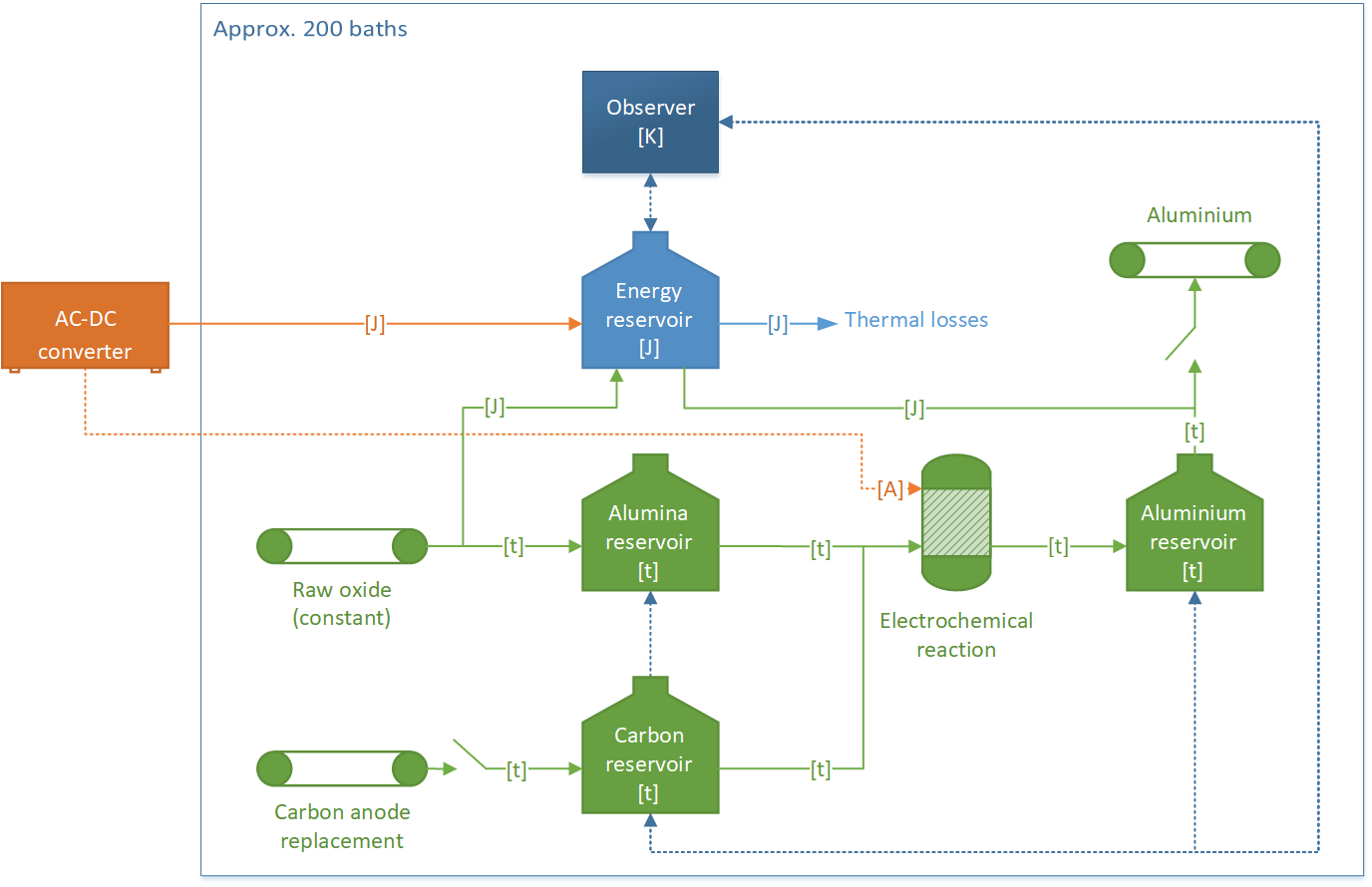}
\par\end{centering}
\caption{\label{fig:Reservoirs-electrolysis}Electrolysis modelled with reservoirs.}
\end{figure}

As a consequence, the electrolytic bath model contains four reservoirs:
\begin{enumerate}
\item The energy is included as a \emph{decaying reservoir}, exactly like
in the previous models \eqref{eq:HVAC-model-state} and \eqref{eq:Oven-model-state-energy},
with some contribution due to the alumina input $\Delta\mathbf{m_{Al_{2}O_{3}}}_{t}^{+}$
and to the metallic aluminium output $\Delta\mathbf{m_{Al}}_{t}^{-}$,
but also to the fact that the reaction is exothermic (with $\Delta H$
being the enthalpy change due to the reaction): 
\begin{align}
\mathbf{E}_{t+1} & =\mathbf{E}_{t}+\mathbf{heating}_{t}-\mathbf{losses}_{t}\label{eq:Electrolysis-model-state-energy}\\
 & +C_{p}^{\mathrm{Al_{2}O_{3}}}\times\Delta\mathbf{m_{Al_{2}O_{3}}}_{t}^{+}\times T_{t}^{\mathrm{out}}\nonumber \\
 & -C_{p}^{\mathrm{Al}}\times\Delta\mathbf{m_{Al}}_{t}^{-}\times\mathbf{T}_{t}\nonumber \\
 & +\Delta H\times\Delta\mathbf{m_{Al}}_{t}^{+}\times\left(\mathbf{T}_{t}-T_{t}^{\mathrm{out}}\right)
\end{align}
Again, the coefficient $C_{p}^{X}$ is the specific heat of material
$X$. The thermal losses can be expressed as linear with respect to
the temperature difference with the outside: 
\[
\mathbf{losses}_{t}=k\times\left(\mathbf{T}_{t}-T_{t}^{\mathrm{out}}\right).
\]
\item The carbon anode \emph{reservoir} can only be consumed by the electrochemical
reaction (it is replaced during maintenance, which we do not mean
to optimise): 
\begin{equation}
\mathbf{m_{C}}_{t+1}=\mathbf{m_{C}}_{t}-\Delta\mathbf{m_{C}}_{t}^{-}.\label{eq:Electrolysis-model-state-carbon}
\end{equation}
\item The alumina \emph{reservoir} is constantly fed, and consumed by the
electrochemical reaction: 
\begin{equation}
\mathbf{m_{Al_{2}O_{3}}}_{t+1}=\mathbf{m_{Al_{2}O_{3}}}_{t}+\Delta\mathbf{m_{Al_{2}O_{3}}}_{t}^{+}-\Delta\mathbf{m_{Al_{2}O_{3}}}_{t}^{-}.\label{eq:Electrolysis-model-state-alumina}
\end{equation}
\item The aluminium \emph{reservoir} is filled by the electrochemical reaction
and periodically emptied: 
\begin{equation}
\mathbf{m_{Al}}_{t+1}=\mathbf{m_{Al}}_{t}+\Delta\mathbf{m_{Al}}_{t}^{+}-\Delta\mathbf{m_{Al}}_{t}^{-}.\label{eq:Electrolysis-model-state-aluminium}
\end{equation}
\item The temperature \emph{observer }depends on the complete mass within
the electrolytic bath. A simplifying assumption is to consider that
the temperature is uniform within the bath: 
\begin{equation}
\mathbf{E}_{t}=\left(C_{p}^{\mathrm{Al_{2}O_{3}}}\times\mathbf{m_{Al_{2}O_{3}}}_{t}+C_{p}^{\mathrm{C}}\times\mathbf{m_{C}}_{t}+C_{p}^{\mathrm{Al}}\times\mathbf{m_{Al}}_{t}\right)\,\mathbf{T}_{t}.\label{eq:Electrolysis-model-temperature}
\end{equation}
As for an oven \eqref{eq:Oven-model-temperature-bounds}, the temperature
bounds can be written linearly based on this expression: 
\[
T_{\min}\times\left(C_{p}^{\mathrm{Al_{2}O_{3}}}\times\mathbf{m_{Al_{2}O_{3}}}_{t}+C_{p}^{\mathrm{C}}\times\mathbf{m_{C}}_{t}+C_{p}^{\mathrm{Al}}\times\mathbf{m_{Al}}_{t}\right)\leq\mathbf{E}_{t},
\]
\[
\mathbf{E}_{t}\leq T_{\max}\times\left(C_{p}^{\mathrm{Al_{2}O_{3}}}\times\mathbf{m_{Al_{2}O_{3}}}_{t}+C_{p}^{\mathrm{C}}\times\mathbf{m_{C}}_{t}+C_{p}^{\mathrm{Al}}\times\mathbf{m_{Al}}_{t}\right).
\]
\item The last block within the process, the electrochemical reactions,
are modelled as \emph{processes}, whose conversion rates $k_{X}^{\prime\prime}$
depend on the average current through the bath ($\mathbf{\overline{I}}_{t}$):
\[
\Delta\mathbf{m_{Al}}_{t}^{+}=k_{\ce{Al}}^{\prime\prime}\times\mathbf{\overline{I}}_{t},\qquad\Delta\mathbf{m_{C}}_{t}^{-}=k_{\ce{C}}^{\prime\prime}\times\mathbf{\overline{I}}_{t},
\]
\[
\Delta\mathbf{m_{Al_{2}O_{3}}}_{t}^{-}=k_{\ce{Al}_{2}\text{O}_{3}}^{\prime\prime}\times\mathbf{\overline{I}}_{t}.
\]
\item The AC-DC converter must be modelled as an \emph{external process}
to provide the required power to heat the bath; moreover, the actual
amperage must be explicitly represented in the model for the electrochemical
reaction rate: 
\[
\mathbf{heating}_{t}=f\!\left(\mathbf{electricity}_{t},\mathbf{\overline{I}}_{t}\right).
\]
\item Finally, the initial conditions are set according to the way the process
is managed: the process is rarely started from scratch (with a temperature
equal to $T_{0}^{\mathrm{out}}$), but rather continuously operated.
\end{enumerate}
The model is shown in Figure~\ref{fig:Reservoirs-electrolysis}.
Existing models tend to use electrochemical multiphysics techniques~\citep{Hofer2011,Zhao2020},
which are not well-suited for light applications like flexibility
estimation.
\begin{remark}
The observer equation \ref{eq:Electrolysis-model-temperature} might
consider directly the sum of all three material reservoirs (alumina,
carbon, and metallic aluminium), without a distinction between the
various materials: their specific heats are of the same order of magnitude
(between $700$ and $900\,\mathrm{J/K\,kg}$). This approximation
is justified in the context of flexibility, as the model does not
need to be very precise. This simplification is useful when fitting
the model to actual data, as fewer parameters must be estimated.
\end{remark}
\begin{remark}
If the electrical current is not considered for flexibilisation, then
the model can be simplified, as the electrochemical reaction then
becomes constant: the mass reservoir levels only change when the bath
is powered on (and thus when the electrochemical reaction takes place).
The only possible values are then naturally discrete.
\end{remark}

\section{Nonlinearity \label{sec:Nonlinearity}}

Another great advantage of the proposed methodology is that most models
are linear. However, these can be limiting for some behaviours that
cannot be completely represented with linear equations. For those
cases, such as ovens (for which a conceptual model is presented in
Section~\ref{subsec:Oven-model}), nonlinearities can be introduced
in the reservoir formalism.

Those nonlinearities can be dealt with in different ways.
\begin{itemize}
\item \textbf{Use of nonlinear nonconvex solvers. }Nonlinear equalities
such as~\eqref{eq:Electrolysis-model-temperature} make standard
optimisation solvers unusable because of their nonconvexity. Therefore,
the solvers must be able to deal with nonconvex constraints and (often)
mixed-integer variables. Programs like Couenne~\citep{Couenne2009},
POD~\citep{Nagarajan2016} (global solvers, not relying on the convexity
assumption), or Bonmin~\citep{Bonmin2013} (\emph{approximating}
the problem as convex) must then be used. The most recent versions
of standard optimisation solvers like CPLEX~\citep{IBM2020} and
Gurobi~\citep{GurobiOptimization2020} now allow for certain types
of nonconvexity (CPLEX 12.7 accepts nonconvex quadratic objective
functions, whereas Gurobi 9.0 also tackles quadratic nonconvex constraints).
\item \textbf{Linear reformulations. }Another option is to formulate the
problem linearly. The main technique we used for this is discretisation:
instead of having continuous variables, some of them are only allowed
to take discrete levels. The choice of variables for discretisation
is made so that, when these variables take a fixed value, the constraint
becomes linear. For example, the constraint~\eqref{eq:Electrolysis-model-temperature}
is linear once the mass is known. Binary variables are then used to
choose among the various discrete values, which keeps the overall
linearity.
\end{itemize}

\subsection{Discretisation \label{subsec:Discretisation-basic}}

Discretisation is an approximation of the previous models. It can
however be used in some cases without being a rough estimate, due
to the actual operational conditions: for instance, industrial ovens
are often loaded by batches, whose sizes can be used to define the
discretisation levels.

An oven like in Section~\ref{subsec:Oven-model} can be reformulated
with this discretisation approach, more specifically the temperature
definition~\eqref{eq:Oven-model-temperature}. Instead of having
a continuous variable $\mathbf{m}_{t}$, it takes its value in a discrete
set $\mathcal{M}=\left\{ \mathcal{M}_{i}\,|\,i\in\left[1,M\right]\right\} $.
For example, instead of having an oven whose mass may freely vary
between $500\,\mathrm{kg}$ and $5\,\mathrm{t}$, its discretised
version may only take the values in $\mathcal{M}=\left\{ 500,1500,2500,3500,5000\right\} \,\mathrm{kg}$,
or $\mathcal{M}=\left\{ 500,1000,1500,2000,3500,5000\right\} \,\mathrm{kg}$
if more precision is needed for low masses.

In order to linearise the quotient $\mathbf{T}_{t}=\mathbf{E}_{t}/C_{p}\,\mathbf{m}_{t}$,
the following expression can be used: 
\begin{align*}
\mathbf{T}_{t} & =\frac{\mathbf{E}_{t}}{C_{p}\,\mathbf{m}_{t}},\qquad\forall t\in\mathcal{T}\\
 & =\begin{cases}
\mathbf{E}_{t}/C_{p}\,\mathcal{M}_{1} & \mbox{ if }\mathbf{m}_{t}=\mathcal{M}_{1}\\
\mathbf{E}_{t}/C_{p}\,\mathcal{M}_{2} & \mbox{ if }\mathbf{m}_{t}=\mathcal{M}_{2}\\
\vdots\\
\mathbf{E}_{t}/C_{p}\,\mathcal{M}_{M} & \mbox{ if }\mathbf{m}_{t}=\mathcal{M}_{\text{\ensuremath{\left|\mathcal{M}\right|}}}
\end{cases},\qquad\forall t\in\mathcal{T}.
\end{align*}

Implementing this in a mathematical optimisation model can be done
in a classical way~\citep{Liberti2009,Liberti2009a,Vielma2015}.
New variables are introduced:
\begin{itemize}
\item The binary variables $\mathbf{m}_{t}^{\left(i\right)}$ indicates
whether the material contents of the oven $\mathbf{m}_{t}$ take the
discrete value $\mathcal{M}_{i}$.
\item The continuous variable $\mathbf{T}_{t}$ is defined as the quotient
$\mathbf{E}_{t}/C_{p}\,\mathbf{m}_{t}$.
\item The continuous variable $\mathbf{T}_{t}^{\left(i\right)}$ is defined
as the quotient $\mathbf{E}_{t}/C_{p}\,\mathbf{m}_{t}$ \emph{if}
$\mathbf{m}_{t}=\mathcal{M}_{i}$, and zero otherwise. These variables
could be called ``partial temperatures'', as their sum yields the
temperature.
\end{itemize}
Then, constraints are added to impose those semantics.
\begin{itemize}
\item Exactly one discrete mass value is possible at any time step: 
\[
\sum_{i\in\mathcal{M}}\mathbf{m}_{t}^{\left(i\right)}=1,\qquad\forall t\in\mathcal{T}.
\]
\item The mass is then defined as a linear combination of those binary choices:
\[
\mathbf{m}_{t}=\sum_{i\in\mathcal{M}}\mathbf{m}_{t}^{\left(i\right)}\,\mathcal{M}_{i},\qquad\forall t\in\mathcal{T}.
\]
\item Similarly, the temperature is given by the sum over all partial temperatures:
\[
\mathbf{T}_{t}=\sum_{i\in\mathcal{M}}\mathbf{T}_{t}^{\left(i\right)},\qquad\forall t\in\mathcal{T}.
\]
\item The partial temperatures $\mathbf{T}_{t}^{\left(i\right)}$ are linked
to the binary choices by the following upper and lower bounds: 
\[
\mathbf{T}_{t}^{\left(i\right)}\geq\frac{\mathbf{E}_{t}}{C_{p}\,\mathcal{M}_{i}}-\frac{E_{\max}}{C_{p}\,m_{\mathrm{min}}}\,\left(1-\mathbf{m}_{t}^{\left(i\right)}\right),\qquad\forall t\in\mathcal{T},\quad\forall i\in\mathcal{M},
\]
\[
\mathbf{T}_{t}^{\left(i\right)}\leq\frac{\mathbf{E}_{t}}{C_{p}\,\mathcal{M}_{i}},\qquad\forall t\in\mathcal{T},\quad\forall i\in\mathcal{M},
\]
\[
\mathbf{T}_{t}^{\left(i\right)}\leq\frac{E_{\max}}{C_{p}\,m_{\min}}\,\mathbf{m}_{t}^{\left(i\right)},\qquad\forall t\in\mathcal{T},\quad\forall i\in\mathcal{M}.
\]
The two last bounds are required to be specified separately, to ensure
that a partial temperature is forced to be zero if the corresponding
mass level is chosen, and takes its expected value otherwise.
\end{itemize}
Thanks to this technique, the temperature definition becomes linear.
It could be extended to more general expressions.

This formulation can be strengthened in order to improve solving times.
The bounds on the partial temperatures $\mathbf{T}_{t}^{\left(i\right)}$
should be as tight as possible to keep good solving times; this issue
motivated the choice of temperature $\mathbf{T}_{t}=E_{t}/C_{p}\,\mathbf{m}_{t}$
over the raw quotient of variables $\mathbf{E}_{t}/\mathbf{m}_{t}$,
as the factor $1/C_{p}$ can reduce the values that are considered
by several orders of magnitude.
\begin{itemize}
\item The bounds on the partial temperatures also have lower bounds, using
the same binary variables: 
\[
\mathbf{T}_{t}^{\left(i\right)}\geq\frac{E_{\mathrm{min}}}{C_{p}\,m_{\mathrm{max}}}\,\mathbf{m}_{t}^{\left(i\right)},\qquad\forall t\in\mathcal{T},\quad\forall i\in\mathcal{M}
\]
\item The mutual exclusion of the $m_{t}^{\left(i\right)}$ can be further
imposed with clique constraints: 
\[
\sum_{i\in\mathcal{V}}\mathbf{m}_{t}^{\left(i\right)}\leq1,\qquad\forall t\in\mathcal{T},\quad\forall\mathcal{V}\subseteq\mathcal{M}:\:\left|\mathcal{V}\right|\geq2
\]
\end{itemize}

\subsection{Benchmark}

The two nonlinear models are benchmarked against each other, in order
to compare their performance when solving the same problem. We use
an oven-like model to heat a given quantity of metal that may be added
at any rate and time (which gives a nonlinear model). Fifty time steps
are considered. Heating happens with temperature ramping constraints.
All models have been written using JuMP~\citep{Dunning2017} in Julia~\citep{Bezanson2017}.

Those models are all compared to a base line, which heats the material
as soon as possible (while respecting the same ramping constraints)
and keeps it at the right temperature until the end of the horizon.

The results are shown in Table~\ref{tab:Nonlinearity-models-comparison}.
When a high precision is needed in the discretised variables, both
CPLEX~\citep{IBM2020} and Gurobi~\citep{GurobiOptimization2020},
two state-of-the-art mixed-integer linear solvers, have troubles to
reach a very low gap, albeit their solutions are more than satisfying
for industrial applications. In the same time budget, the nonlinear
formulation with open-source nonlinear solvers achieves a better solution,
without being hindered by discretised variables.

\begin{table}
\caption{\label{tab:Nonlinearity-models-comparison}Solving time for the various
nonlinear models. All solvers were stopped after approximately one
hour of computations, even when they did not find the optimum solution
(nonzero gap). Those tests were run on a machine with two Intel Xeon
E5-2650v4 (2.2GHz) and 128GB of RAM.}

\begin{tabular}{>{\raggedright}p{3cm}lll>{\raggedright}p{2cm}>{\raggedright}p{2cm}}
\hline 
Model & Underlying solver & Time (s) & Solution cost (€) & Gap reported by the solver (\%) & Cost improvement with respect to reference (\%)\tabularnewline
 &  &  &  &  & \tabularnewline
\hline 
Reference: reach target temperature as soon as possible & CPLEX 12.7.1 & $0.02$ & $10,793.25$ & $0.00$ & \textemdash{}\tabularnewline
\hline 
\multirow{5}{3cm}{Reservoir model: linearisation by discretisation (16 levels)} & CPLEX 12.7.1 & $3,610.37$ & $8,208.23$ & $2.68$ & $23.95$\tabularnewline
 & Gurobi 7.5.0 & $3,600.02$ & $8,034.72$ & $2.12$ & $25.56$\tabularnewline
 &  &  &  &  & \tabularnewline
 &  &  &  &  & \tabularnewline
 &  &  &  &  & \tabularnewline
\hline 
\multirow{4}{3cm}{Reservoir model: nonlinear (continuous)} & Bonmin 1.8.4 & $3,600.59$ & $7,259.46$ & $\text{\textemdash}$ & $32.74$\tabularnewline
 & Couenne 0.5.4 & $3,601.39$ & $7,260.20$ & $100$ & $32.73$\tabularnewline
 &  &  &  &  & \tabularnewline
 &  &  &  &  & \tabularnewline
\hline 
\end{tabular}
\end{table}

\section{Process typology \label{sec:Process-typology}}

Based on these models, we can derive a typology of industrial processes
based on a few characteristics. Table~\ref{tab:Typology-industrial-processes}
does so according to two criteria, focusing on the \emph{flows of
material to process}:
\begin{itemize}
\item the \emph{inputs} to the process: are they continuous or periodical?
\item the \emph{outputs} from the process: are they continuous or periodical?
\end{itemize}
Those two questions can be equivalently formulated as follows. Is
the mass currently processed constant? How can it vary (continuously,
periodically)?

Combined, these three parameters indicate how the process could be
modelled. A mass that is constant (HVAC) or highly predictable (electrolysis,
batches; kilns to a lesser degree) can result in model simplifications.
Completely variable masses often implies nonlinearity (as most ovens,
see Section~\ref{sec:Nonlinearity}).

\begin{table}
\caption{\label{tab:Typology-industrial-processes}Typology of industrial processes
to help reservoir modelling.}

\begin{tabular}{>{\raggedright}p{4cm}|l>{\raggedright}p{2.5cm}>{\raggedright}p{4cm}}
\hline 
Process & Mass & Inputs & Outputs\tabularnewline
 &  &  & \tabularnewline
\hline 
HVAC & Constant & \multicolumn{2}{l}{(Constant mass: no inputs nor outputs)}\tabularnewline
\hline 
Kiln & Variable & Continuous & Continuous\tabularnewline
\hline 
Oven (including electric arc furnace, induction furnace) & Batch & Discrete (batches: when starting) & Discrete (batches: when done, e.g. molten metal)\tabularnewline
\hline 
Electrolysis & Variable & Continuous & Discrete (periodical: e.g., every day)\tabularnewline
\hline 
\end{tabular}
\end{table}

The diagrams shown in Section~\ref{sec:Reservoir-models} also help
build another typology, this time based on the \emph{potential flexibility}
levers for each process. Indeed, the presence of a decaying heat energy
reservoir indicates that the process is amenable to load shifting:
even if heating is stopped, the process might still continue to produce,
albeit probably at a lessened rate. If multiple fuels are possible
for heating, then fuel switching can be used. For all non-continuous
processes, load scheduling can be applied; for continuous processes
whose production can be tuned (like electrolysis), exploiting flexibility
may lead to load shedding. The processes that have been studied in
Section~\ref{sec:Reservoir-models} are included in Table~\ref{tab:Analysis-flexibility-levers}.

\begin{table}
\caption{\label{tab:Analysis-flexibility-levers}Analysis of flexibility levers
available in industrial processes.}

\begin{tabular}{l|>{\raggedright}p{8cm}}
\hline 
Process & Potential flexibility levers\tabularnewline
 & \tabularnewline
\hline 
HVAC & Thermal inertia: \emph{load shifting}\tabularnewline
Kiln & Thermal inertia: \emph{load shifting} (production delayed)\\
Starting time: \emph{load scheduling}\\
\emph{Fuel switching}\tabularnewline
Continuous oven & Thermal inertia: \emph{load shifting} (production delayed)\emph{}\\
\emph{Fuel switching}\tabularnewline
EAF & Thermal inertia: \emph{load shifting} (production delayed)\\
Starting time: \emph{load scheduling}\tabularnewline
Electrolysis & Thermal inertia: \emph{load shifting}\\
For a long period of time, becomes \emph{load shedding}!\tabularnewline
 & \tabularnewline
\hline 
\end{tabular}
\end{table}

\section{Flexibility potential of industrial processes \label{sec:Results}}

Some models developed in Section~\ref{sec:Reservoir-models} are
now fit to industrial data, and the flexibility potential of the processes
is estimated based on historical scenarios. A major hypothesis is
that exploiting the flexibility for these industrial sites has no
impact on the electricity market, as they do not consume enough electricity.

All optimisation models have been written using JuMP~\citep{Dunning2017}
in Julia~\citep{Bezanson2017}. The source code for the simulations
is available online at the following address:\\
\url{https://github.com/dourouc05/IndustrialProcessFlexibilisation.jl}

\subsection{Induction furnace \label{subsec:Results-Induction-furnace}}

The induction furnace of Section~\ref{subsec:Induction-furnace}
needs three parameters. The first one, $\alpha$, characterises the
energy level depending on the initial temperature of the material.
The second one, $\beta$, is the efficiency of converting the electrical
power into heat (see Appendix \ref{sec:Appendix-IF}). The last one,
$\delta$, indicates the energy losses.

Each sample $s\in\mathcal{S}$, taken from historical measurements,
corresponds to one use of the furnace. It mainly consists in a power
curve $P_{s}^{t}$, indicating the average power injected through
the circuits each hour (a melting cycle lasts twelve hours). To have
a better fit, each sample may have its own value of $\beta$, denoted
by $\beta_{s}$. The values of the $\beta_{s}$ are brought closer
together by a constraint limiting the variance of the $\beta_{s}$
to $0.001$. All in all, fitting the parameters is done through the
following optimisation program (a convex QCQP), based on the reservoir
model of Section~\ref{subsec:Induction-furnace}: 
\[
\begin{array}{ccc}
\min & \sum_{s\in\mathcal{S}}\left(\mathbf{E}_{s}^{t}-E_{s}\right)^{2}\\
\text{s.t.} & \mathbf{E}_{s}^{0}=\boldsymbol{\alpha}\times m_{s}\times T_{0} & \forall s\in\mathcal{S}\\
 & \mathbf{E}_{s}^{t}=\mathbf{E}_{s}^{t-1}+\boldsymbol{\beta_{s}}\times P_{s}^{t}-\boldsymbol{\delta} & \forall s\in\mathcal{S},\forall t\in\mathcal{T}\\
 & \left|\mathcal{S}\right|\times\boldsymbol{\overline{\beta}}=\sum_{s\in\mathcal{S}}\boldsymbol{\beta_{s}}\\
 & \sum_{s\in\mathcal{S}}\left(\boldsymbol{\beta_{s}}-\mathbf{\boldsymbol{\overline{\beta}}}\right)^{2}\leq0.001\\
 & \boldsymbol{\alpha}\geq0\\
 & \boldsymbol{\beta}_{s}\geq0 & \forall s\in\mathcal{S}\\
 & \mathbf{\boldsymbol{\overline{\beta}}}\geq0\\
 & \boldsymbol{\delta}\geq0
\end{array}
\]
To evaluate the results of this model, we use a leave-one-out procedure
\citep{Efron1982}: for each sample $s\in\mathcal{S}$, we solve the
previous optimisation program over the samples $\mathcal{S}\backslash\left\{ s\right\} $
to get a value $\beta_{\mathcal{S}\backslash\left\{ s\right\} }$,
and we average the obtained square errors to predict the energy of
$s$ based on $\beta_{\mathcal{S}\backslash\left\{ s\right\} }$.
It results in a root mean square error on our data set of $65.194\,\text{kWh}$
(the ratio of this error to the average energy is $1.78\%$).

Once a reservoir model is fit, it can be used to estimate the flexibility
potential of the process. To this end, we compare our reservoir-based
optimisation to current fixed consumption profiles on a historical
price scenario (an average day of January 2016 on the Belgian day-ahead
market). No optimisation is currently performed on the schedule: the
smelting process always starts at the same hour and lasts for twelve
hours.

Related constraints must be added to the formulation in order to implement
real-world constraints. Mostly, the peak power is limited, and the
heating cannot abruptly change. Three constraints are thus added for
each time step: a minimum and a maximum power ($P_{\min}$ and $P_{\max}$,
respectively), and a ramping constraint ($\rho_{\min}$ and $\rho_{\max}$
are respectively the minimum and maximum ratios by which the electrical
consumption is allowed to change from hour to hour).

The complete model is therefore the following, where $p_{h}$ indicates
the price of electricity for the hour $h$ (in €/MWh), $c_{h}$ the
electrical consumption (in MWh), and $E_{h}$ the energy of the metal
to melt: 
\[
\begin{array}{ccc}
\min & \sum_{h=1}^{H}p_{h}\,c_{h}\\
\text{s.t.} & E_{0}=\alpha\,m\,T_{0},\\
 & E_{H}=\alpha\,m\,T_{H},\\
 & E_{h}=E_{h-1}+\beta\,c_{h}-\delta & \forall h\in\left\{ 1,2\dots H-1\right\} \\
 & P_{\min}\leq c_{h}\leq P_{\max} & \forall h\in\left\{ 1,2\dots H\right\} \\
 & \rho_{\min}\,c_{h}\leq c_{h+1}\leq\rho_{\max}\,c_{h} & \forall h\in\left\{ 1,2\dots H\right\} 
\end{array}
\]

\begin{figure}
\begin{centering}
\includegraphics[scale=0.6]{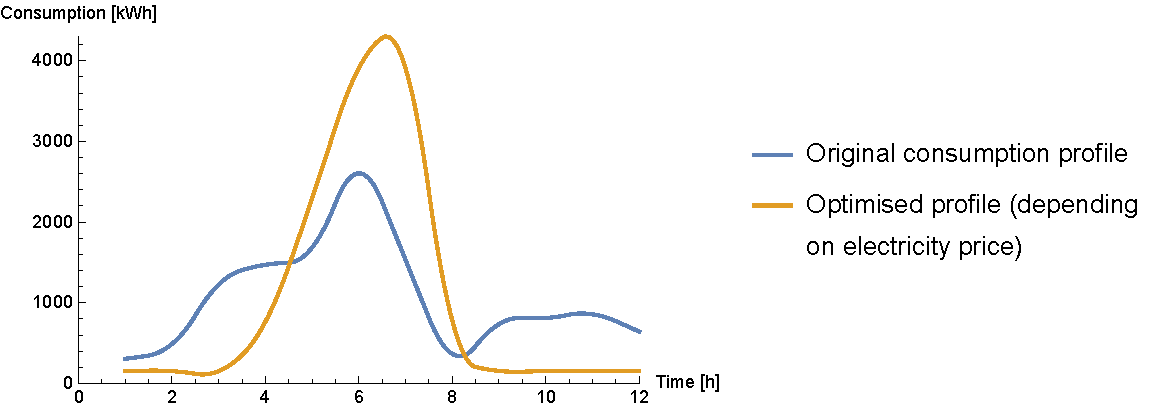}
\par\end{centering}
\caption{\label{fig:Comparison-consumption-profiles}Comparison of consumption
profiles between the existing fixed solution and the result of the
optimisation.}
\end{figure}

Per heat, using a reservoir model to decide the heating power, hour
per hour, could decrease the costs by 8.35\% per heat, from €239.11
to €219.15; over a year, this corresponds to more than €15,000 of
savings. The main difference in the planned power consumption is that
its peak is shifted to exploit the lowest prices during the production
period (as shown in Figure \ref{fig:Comparison-consumption-profiles}).

Computationally speaking, this model allows to optimise the required
power to a given price scenario in a fraction of a second: a 95\%
confidence interval is $0.03\pm0.01$ seconds (with either CPLEX 12.7.1
or Gurobi 7.5.0). Fitting the parameters is also very quick, as the
leave-one-out validation phase takes less than two minutes for the
available data. The from-scratch implementation in Julia takes 150
lines, but the process could be automatised with a GUI to build the
reservoir model. These characteristics are very appealing for energy-sector
consultants, who may want to estimate the flexibility potential of
an industrial site in very little time: most of the effort can be
spent on actual discussions on the flexibility solutions that can
be implemented.

\subsection{Industrial cooling \label{subsec:Results-Industrial-cooling}}

Based on industrial data from a polypropylene-film plant, a model
similar to that of Section~\ref{subsec:Industrial-cooling} can be
tuned. All the needed physical constants are known and match the data
set. The model is then fed with scenarios of water volumes to cool
down. In practice, the price scenarios can be estimated in advance,
while the heat and water volumes are known with a high precision.

The processes in Figure \ref{fig:Reservoirs-Cooling} have to be specified
in more details. A linear model is deemed sufficient for our needs,
using coefficients of performance, and the functions $f$ and $g$
can be written as: 
\[
f\!\left(\mathbf{P}_{\mathrm{CT},t}\right)=\mathrm{cop}_{\mathrm{CT}}\times\mathbf{P}_{\mathrm{CT},t},
\]
\[
f\!\left(\mathbf{P}_{\mathrm{chiller},t}\right)=\mathrm{cop}_{\mathrm{chiller}}\times\mathbf{P}_{\mathrm{chiller},t}.
\]
Thus, the two energy budget constraints \eqref{eq:IC-model-budget-energy-CT}
and \eqref{eq:IC-model-budget-energy-chiller} become: 
\begin{align*}
\mathbf{E}_{\mathrm{CT},t+1} & =\mathbf{E}_{\mathrm{CT},t}+C_{\mathrm{water}}\times T_{\mathrm{p},t}\times\mathbf{water_{p\to CT,t}}\\
 & -\mathbf{E}_{\mathrm{CT\to B},t}-\mathrm{cop}_{\mathrm{CT}}\times\mathbf{P}_{\mathrm{CT},t}.
\end{align*}
\begin{align*}
\mathbf{E}_{B,t+1} & =\mathbf{E}_{B,t}+C_{\mathrm{water}}\times T_{\mathrm{p},t}\times\mathbf{water_{p\to B,t}}+\mathbf{E}_{\mathrm{CT\to B},t}\\
 & -C_{\mathrm{water}}\times\mathbf{T}_{\mathrm{B},t}\times\mathbf{water_{B\to p,t}}-\mathrm{cop}_{\mathrm{chiller}}\times\mathbf{P}_{\mathrm{chiller},t}.
\end{align*}
Two other constraints must be added for the cooling processes, as
they have a limited maximum power: 
\[
\mathbf{P}_{\mathrm{CT},t}\leq P_{\mathrm{CT},\max},
\]
\[
\mathbf{P}_{\mathrm{chiller},t}\leq P_{\mathrm{chiller},\max}.
\]

\begin{remark}
Variable coefficients of performance are not considered in this case
in order to keep the model simple. The increased complexity is not
justified in a context of \emph{approximate} models.
\end{remark}
In order to compare the impact of flexibility on the cooling behaviour,
the temperature bounds are set in two different ways:
\begin{itemize}
\item either $T_{\mathrm{target}}\pm0.5\text{°C}$ (\emph{low-flexibility
scenario}), as currently implemented in the studied use case
\item or $T_{\mathrm{target}}\pm3\text{°C}$ (\emph{high-flexibility scenario}),
the maximum temperature variation that the production equipment may
tolerate, based on a deeper analysis of their data sheets
\end{itemize}
The model includes three nonlinear constraints: \eqref{eq:IC-model-concordance-CT},
\eqref{eq:IC-model-concordance-B}, and \eqref{eq:IC-model-concordance-flow}.
They are implemented as nonlinear constraints; we also compare this
formulation to a finely-discretised linear version. In practice, the
nonconvex formulation can be solved to optimality faster than the
discretised version with a time horizon of eight hours (Table~\eqref{tab:Results-IC-Comparison-formulations}).

In order to compare these two flexibility scenarios, a rolling-horizon
algorithm is implemented. This choice helps keep the running times
low: the optimisation program is run for eight hours, then the result
for the first time step is used as the initial condition for the next
program, whose horizon is shifted by one time step. Even though it
is closer to the plant operating conditions, and therefore better
estimates the actual flexibility potential, it no more guarantees
a global optimality over the complete time horizon.

\begin{table}
\caption{\label{tab:Results-IC-Comparison-formulations}Comparison between
the two formulations, for a time horizon of eight hours, with a time
step of one hour. Twenty-five discretisation steps are used for the
temperature. The programs are solved within a rolling-horizon algorithm;
averages are based on its iteration. The mixed-integer linear model
could not be solved over the one-week horizon due to memory problems
(the formulation is too large). Those tests were run with CPLEX 12.7.1
(MILP formulation), Couenne 0.5.4, and Gurobi 9.0.0 (nonconvex formulation),
on a machine with two Intel Xeon E5-2650 (2.2GHz) and 128GB of RAM.}

\begin{tabular}{lccc}
\hline 
 & Mixed-integer linear & \multicolumn{2}{c}{Nonconvex}\tabularnewline
 & CPLEX 12.7.1 & Couenne 0.5.4 & Gurobi 9.0.0\tabularnewline
\hline 
Number of constraints & 1120 & \multicolumn{2}{c}{100}\tabularnewline
Number of variables & 480 (378 integers) & \multicolumn{2}{c}{102 (0 integers)}\tabularnewline
Average number of explored nodes & 1,193,050,714 & 2,210,153 & 662,159\tabularnewline
Average solving time (seconds) & 62,742.96 & 3,694.52 & 46.95\tabularnewline
 &  &  & \tabularnewline
\hline 
\end{tabular}
\end{table}

\begin{figure}
\begin{centering}
\includegraphics[scale=0.3]{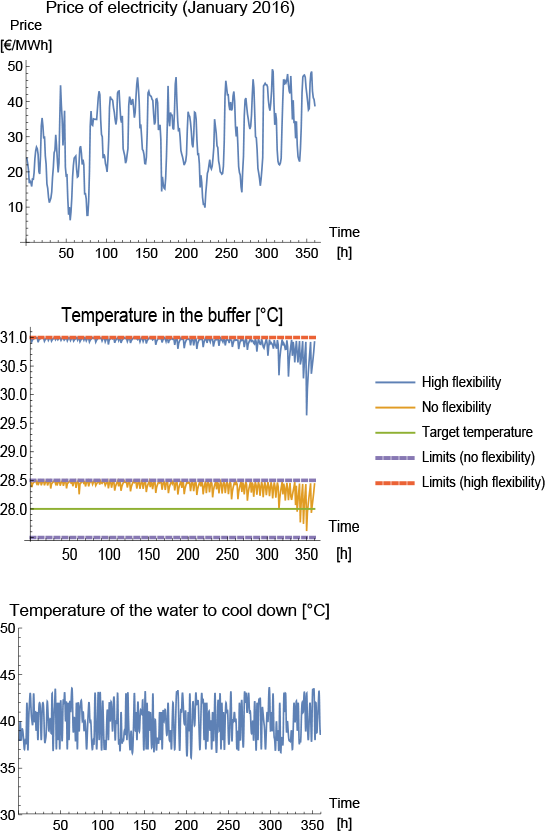}
\par\end{centering}
\caption{\label{fig:Cooling-scenario}Obtained experimental result for one
scenario: increasing flexibility lowers the energy costs by 24\%.}
\end{figure}

When using this methodology on a synthetic heat inflow to cool down\footnote{No historical scenario could be shown in article due to the data being
proprietary.}, the energy costs can be lowered by 24\% when going from the low-
to the high-flexibility scenario (Figure~\ref{fig:Cooling-scenario}):
it goes from €123,569 down to €93,909 for one week (with a price scenario
corresponding to the first week of January 2016 on the Belgian day-ahead
market). Results on historical data are highly similar. What is more,
the obtained solution uses the extra flexibility to lower the temperature
before an increase in both the electricity price and the heat to eliminate:
this is exactly the expected kind of solution.

Similarly to the induction furnace (Section~\eqref{subsec:Results-Induction-furnace}),
computation times are very encouraging: a complex, real-sized, nonconvex
model can be solved quickly to optimality (albeit only using Gurobi
9.0's new nonconvex functionalities). Model-building times are again
reduced to a very low amount: each model corresponds to 100 lines
of Julia when implemented from scratch, with very little thought required
to build the reservoir model, even for someone who is not a specialist
of these processes.

\section{Conclusion \label{sec:Conclusion}}

The reservoir framework helps build models of actual industrial processes.
It does so by following a physically-based approach, around the concept
of storage (either material or energy). This principle corresponds
to the industrial and physical reality, as these models can correctly
approximate the behaviour of processes. Such reservoir models can
lend themselves to interpretation about flexibility (as performed
in Section~\ref{sec:Process-typology}).

This framework defines a series of blocks that can be used to think
about processes. However, it does not specify anything about the actual
energy consumption. For now, they are supposed to be easily represented
by so-called processes. They may hide different kinds of mathematical
expressions, from a constant efficiency to translate fuel consumption
into energy (which is probably the best fit for this framework, due
to its simplicity) to more complex models.

Simple extensions to the proposed building blocks allows modelling
more processes. For example, a crusher can correspond to a reservoir
that transforms one product (such as rock) into another one (like
gravel), with a conversion rate that depends on the energy consumption.
The major problem of this kind of formulation is that linearity is
forgone. More work is needed to look into these kinds of generalisations,
and the way to produce linear models from slightly different building
blocks.

Nevertheless, without needing many alterations, the reservoirs can
already be used to study the impact of flexibility on some industrial
processes, as done in Section~\ref{sec:Results}. More generally,
reservoir-based models share many similarities and are still able
to deal with many different processes, while being lightweight.

\section{Appendix: mathematical derivation of the induction furnace's electrical
efficiency\label{sec:Appendix-IF}}

The power dissipated by eddy currents can be expressed as a linear
function of the electrical current in the coil. Indeed, the dissipation
is related to the maximum magnetic field $B_{\mathrm{max}}$ (if it
varies sinusoidally in time)~\citep{Serway2010}:

\texttt{
\[
P_{\mathrm{eddy}}=\frac{\pi^{2}}{6}\,\frac{e^{2}\,B_{\mathrm{max}}^{2}\,f^{2}}{\rho},
\]
}where $e$ is thickness, $f$ current frequency, $\rho$ material
resistivity. In an induction furnace, $f$ is fixed, $\rho$ depends
on material, the other parameters are a function of the furnace's
geometry. The only variable is thus the maximum magnetic field, which
can be given by Ampère's law, as the coil corresponds to a solenoid
with $N$ wire turns and a total wire length of $\ell$: 
\[
B=\mu\,\frac{N}{\ell}\,I.
\]
$\mu$ is the permeability of the metal to heat and the bucket, and
$I$ is the intensity of the electrical current. This formula can
be rewritten to exhibit the electrical power $P$ in the coil instead
of the current $I$, using Ohm's law $P=R\,I^{2}$: 
\[
B=\mu\,\frac{N}{\ell}\,\frac{\sqrt{P}}{R}.
\]
Finally, the dissipated power in the metal is linked to the injected
electrical power by:

\texttt{
\begin{align*}
P_{\mathrm{eddy}} & =\frac{\pi^{2}}{6}\,\frac{e^{2}\,f^{2}}{\rho}\,\left(\mu\,\frac{N}{\ell}\,\frac{\sqrt{P}}{R}\right)^{2}\\
 & =\frac{\pi^{2}}{6}\,\frac{e^{2}\,f^{2}}{\rho}\,\mu^{2}\,\frac{N^{2}}{\ell^{2}}\,\frac{P}{R^{2}}\\
 & =\underbrace{\frac{\pi^{2}\,e^{2}\,f^{2}\,\mu^{2}\,N^{2}}{6\,\rho\,\ell^{2}\,R^{2}}}_{\mathrm{constant}}\,P.
\end{align*}
}In other words, the power dissipated in the metal is directly proportional
to the electrical power fed into the circuit, with a constant ratio.
\begin{acknowledgements}
We would like to thank Mr. Nicolas Descouvemont (N-SIDE at the time
of the collaboration, now with OMP), Dr. Olivier Devolder (N-SIDE),
and Pr. Quentin Louveaux (university of Liège) for their help in building
and assessing the accuracy of our models.

\noindent\textbf{Conflicts of interest$\quad$}The author declare
that there is no conflict of interest.
\end{acknowledgements}

\bibliographystyle{spbasic}
\addcontentsline{toc}{section}{\refname}\bibliography{C:/Users/Thibaut/Documents/library}

\end{document}